
\magnification=1200
\baselineskip 18pt
\def\I{I\!\!P}
\def\medtype{\let\rm=\tenrm
\baselineskip=12pt minus 1pt
\rm}
\font\ms=msbm10
\font\msp=msbm5
\font\msq=msbm7

\def\Z{\ms Z}

\def\qed{\vrule height 1.1ex width 1.0ex depth -.1ex }
\def\mapright#1{\smash{\mathop{\longrightarrow}\limits^{#1}}}

\centerline{\bf ON THE CORANK OF GAUSSIAN MAPS}
\centerline{\bf FOR GENERAL EMBEDDED K3 SURFACES}
\vskip .5cm
\centerline{Ciro Ciliberto*\footnote\null{{\it \noindent 1991
Mathematics Subject Classification:} Primary 14J28. Secondary 14J10,
14C34.}, Angelo
Felice\footnote\null{{\it Key words and phrases:} K3 surfaces,
Gaussian
maps.} Lopez\footnote*{Research  partially supported by the MURST
national
project ``Geometria Algebrica";} and Rick\footnote\null{the authors
are
members of GNSAGA of CNR.} Miranda*\footnote*{* Research  partially
supported
by NSF under grant DMS-9403392.}}   \vskip 1cm \midinsert \narrower
\hsize=13.3cm
\medtype{

\noindent ABSTRACT: Let $S_{g}$ be a general prime K3 surface in
$\I^g$ of
genus $g \geq 3$ or a general double cover of $\I^2$ ramified
along a sextic curve for $g = 2$ and $S_{i,g}$ its {\it i}-th
Veronese
embedding. In this
article we
compute the corank of the Gaussian map
$\matrix{\Phi_{{\cal O}_{S_{i,g}}(1)} :
\bigwedge^2 H^0(S_{i,g},{\cal O}_{S_{i,g}}(1))  \to
H^0(S_{i,g},\Omega_{S_{i,g}}^1(2)) \cr}$ for $i \geq 2, g \geq 2$
and $i=1, g
\geq 17$. The main idea is to reduce the surjectivity of
$\Phi_{{\cal O}_{S_{i,g}}(1)}$ to an application of the
Kawamata-Viehweg
vanishing theorem on the blow-up of $S_{i,g} \times S_{i,g}$
along
the
diagonal. This is seen to apply once the hyperplane divisor of
the K3
surface $S_{i,g}$ can be decomposed as a sum of three suitable
birationally
ample divisors. We show that such a decomposition exists when
$i \geq 3$ or on
some K3 surfaces, constructed using the surjectivity of the
period
mapping,
when $i = 1, g \geq 17$ or $i=2, g \geq 7$.}
\endinsert
\vskip 1cm
\hsize=15cm
\magnification=1200
\baselineskip 18pt
\hoffset=0cm

\noindent {\bf 1. INTRODUCTION}
\vskip .5cm
Let $\cal H$ be the disjoint union of the Hilbert schemes of smooth
K3
surfaces in $\I^g$, for all $g \geq  3$. It follows by the
transcendental
theory that
$\cal H$ has some main components ${\cal H}_{g}$ whose general
element
is a
genus
$g$ K3 surface
$S_{g} \subset \I^g$ of degree $2g - 2$ with Picard group generated
by its
hyperplane class (the so called {\it prime} K3 surfaces), while all
the
other components ${\cal H}_{i,g}$ for $i \geq 2$ are obtained by
re-embedding prime K3 surfaces via the {\it i}-th Veronese map.

In this paper we will compute the corank of the Gaussian map
$$\matrix{\Phi_{{\cal O}_{S_{i,g}}(1)} : \bigwedge^2
H^0(S_{i,g},{\cal
O}_{S_{i,g}}(1))  \to H^0(S_{i,g},\Omega_{S_{i,g}}^1(2)) \cr}$$

\noindent on a general embedded K3 surface, that is a K3 surface
$S_{i,g}$ representing a general point of ${\cal H}_{i,g}$ (where we
set
${\cal H}_{1,g} = {\cal H}_{g},  S_{1,g} = S_{g}$ and ${\cal H}_{2}$
is the
family of genus $2$ K3  surfaces, that is double covers of $\I^2$
ramified
along a sextic curve). This is of course equivalent to computing the
corank of the Gaussian map
$\matrix{\Phi_{{\cal O}_{S_{g}}(i)} : \bigwedge^2 H^0(S_{g},{\cal
O}_{S_{g}}(i))  \to H^0(S_{g},\Omega_{S_{g}}^1(2i)) \cr}$.
The main technique that we will use to study this map was
suggested to us by L. Ein and it consists, as we will see in
section 2, in
the fact that its surjectivity follows by the Kawamata-Viehweg
vanishing
theorem once the hyperplane divisor of the K3 surface $S_{i,g}$ can
be
decomposed as a sum of three suitable birationally ample  divisors
(Lemmas
(2.1) and (2.2)). Therefore this approach is particularly effective,
and in
fact it gives sharp results, when such a decomposition exists, that
is on
some surfaces representing points in ${\cal H}_{g}$ for high enough
genus
or in ${\cal H}_{i,g}, i \geq 2$ where the hyperplane divisor is
divisible.
The existence of K3 surfaces whose {\it i}-th Veronese embedding has
hyperplane  divisor decomposable as above will then be treated in
section
3, mainly using the surjectivity of the period mapping.

Our result is as follows:
\vskip .3cm
{\bf Theorem (1.1).} {\sl For $i \geq 1$ and $g \geq 2$ let
${\cal H}_{i,g}$ be the component of the Hilbert scheme of K3
surfaces
whose general element $S_{i,g}$ is the {\it i}-th Veronese
embedding of a
prime genus $g$ K3 surface. Then the corank of $\Phi_{{\cal
O}_{S_{i,g}}(1)}$ is given by the following table:}
\vskip .3cm

\centerline{\vbox{\offinterlineskip
\hrule
\halign{&\vrule#&
\strut\hfil#\tabskip=1.5em\cr
height4pt&\omit&&\omit&&\omit&\cr
&\ \ \ $i$\hfil&&\ \ \ \ $g$\hfil&&\ \ $corank \ \Phi_{{\cal
O}_{S_{i,g}}(1)}$&\cr
height4pt&\omit&&\omit&&\omit&\cr
\noalign{\hrule}
height4pt&\omit&&\omit&&\omit&\cr
&\ \ \ 1\hfil&&\ \ $\geq 17$\hfil&&0\hfil&\cr
height4pt&\omit&&\omit&&\omit&\cr
&\ \ \ 2\hfil&&\ \ \ \ 2\hfil&&3\hfil&\cr
height4pt&\omit&&\omit&&\omit&\cr
&\ \ \ 2\hfil&&\ $\geq 3$\hfil&&0\hfil&\cr
height4pt&\omit&&\omit&&\omit&\cr
&\ \ \ 3\hfil&&\ \ \ \ 2\hfil&&8\hfil&\cr
height4pt&\omit&&\omit&&\omit&\cr
&\ \ \ 3\hfil&&\ $\geq 3$\hfil&&0\hfil&\cr
height4pt&\omit&&\omit&&\omit&\cr
&\ $\geq 4$\hfil&&\ $\geq 2$\hfil&&0\hfil&\cr
} \hrule} \ .}

\noindent {\sl Moreover for all $i,g$ in the above table (with $i
\geq 5$ if
$g = 2$), the Gaussian maps $\Phi_{{\cal O}_{S_{i,g}}(1), {\cal
O}_{S_{i,g}}(k)}$ are surjective for all $k \geq 2$.}

\vskip .5cm
Remark that it remains to compute the corank for $i = 1, g \leq
16$; on
the other hand the above result for $i = 1, g \geq 17$ allows
a generalization of a theorem of [CLM] (see below).

As is well-known now
Gaussian maps give a new interesting approach in  the study of
algebraic
varieties (see for example [W1], [W3], [CLM], [Z]). In particular,
in the
above case, the knowledge of the corank of
$\Phi_{{\cal O}_{S_{i,g}}(1)}$ gives, upon restricting to a smooth
hyperplane section $C_{i,g}$ of $S_{i,g}$, the possibility of
computing
the corank of the Gaussian map (or Wahl map)
$$\matrix{\Phi_{\omega_{C_{i,g}}} : \bigwedge^2 H^0(C_{i,g},
\omega_{C_{i,g}}) \to H^0(C_{i,g},\omega_{C_{i,g}}^{\otimes 3})
\cr}.$$

\noindent The latter is a particularly interesting invariant as it
encodes
information both on the component ${\cal H}_{i,g}$ of the Hilbert
scheme
(since it gives the dimension of the tangent space at points
representing
cones over $C_{i,g}$) and on the possibility of extending $C_{i,g}$
as a
curve section of higher dimensional varieties (via Zak's theorem
[Z],
[BEL]). The main application that we have in mind is to use the
knowledge
of $corank \ \Phi_{\omega_{C_{i,g}}}$ to classify Fano threefolds
of index
greater than one, in the same vein as [CLM] where we carried out
this
project for prime Fano threefolds. The calculation of $corank \
\Phi_{\omega_{C_{i,g}}}$ (that gives an extension, for $g \geq 17$,
of Theorem 4.1 of [CLM]) and the study of its consequences on Fano
threefolds of index greater than one will appear in a forthcoming
paper.
\vskip .5cm
\noindent {\bf 2. EIN'S APPROACH TO THE SURJECTIVITY OF GAUSSIAN
MAPS}
\vskip .5cm
Let $S$ be a smooth embedded K3 surface, ${\cal O}_{S}(1)$ its
hyperplane
bundle. In this section we will give some sufficient conditions for
the
surjectivity of the Gaussian map $\matrix{{\Phi_{{\cal O}_{S}(1)}} :
\bigwedge^2 H^0(S,{\cal O}_{S}(1)) \to H^0(S,\Omega_{S}^1(2)) \cr}$,
which
is defined by
$\Phi_{{\cal O}_{S}(1)} (\sigma \wedge \tau) = \sigma d\tau - \tau
d\sigma$ (on any open subset where ${\cal O}_{S}(1)$ is trivial).

The main idea, that was suggested to us by L. Ein, is to blow up
$S \times S$ along its diagonal $\Delta$ and then use the
Kawamata-Viehweg
vanishing theorem. To this end let us recall that a line bundle $A$
on a
projective variety $X$ is said to be {\it nef} if $c_{1}(A) \cdot
\Gamma
\geq 0$ for every irreducible curve $\Gamma \subset X$; $A$ is
{\it big}
if for some $m > 0$ the rational map defined by $mA$ on $X$ is
birational.
The Kawamata-Viehweg vanishing theorem asserts that $H^i(X,
\omega_{X}
\otimes A) = 0$ for $i > 0$ if $X$ is smooth and
$A$ is big and nef.

Now let $Y$ be the blow-up of $S \times S$ along its diagonal
$\Delta$,
$E$ the exceptional divisor and for every sheaf $\cal F$ on $S$
let us
denote by  ${\cal F}_{i}, i=1,2$, its pull-back via the map
$Y \to S \times S \mapright{p_{i}} S$ where $p_{i}$ is the
{\it i}-th
projection. Then we have

{\bf Lemma (2.1).} {\sl Suppose there are line bundles
$A_{1}, A_{2}, A_{3}$ on $S$ such that

(i) ${\cal O}_{S}(1) \cong A_{1} \otimes A_{2} \otimes A_{3}$;

(ii) $\bigotimes\limits_{j=1}^3 [A_{j1} \otimes A_{j2} (-E)]$ is
big and
nef on $Y$.

\noindent Then $\Phi_{{\cal O}_{S}(1)}$ is surjective.}

\noindent {\it Proof:} Let ${\cal I}_{\Delta}$ be the ideal sheaf
of
$\Delta \subset S \times S$, $L = {\cal O}_{S}(1)$ and consider the
exact
sequence
$$0 \to p_{1}^*L \otimes p_{2}^*L \otimes {\cal I}_{\Delta}^2 \to
p_{1}^*L \otimes p_{2}^*L \otimes {\cal I}_{\Delta} \to
p_{1}^*L \otimes p_{2}^*L \otimes {{\cal I}_{\Delta} \over {\cal
I}_{\Delta}^2} \to 0.$$

\noindent It is a standard fact that the
isomorphism $H^0(\Delta,p_{1}^*L \otimes p_{2}^*L \otimes
{{\cal I}_{\Delta} \over {\cal I}_{\Delta}^2}) \cong
H^0(S,\Omega_{S}^1(2))$ gives
$$Coker \Phi_{L} \cong Coker \{ H^0(S \times S,p_{1}^*L \otimes
p_{2}^*L \otimes
{\cal I}_{\Delta}) \to H^0(\Delta,p_{1}^*L \otimes p_{2}^*L \otimes
{{\cal I}_{\Delta} \over {\cal I}_{\Delta}^2}) \}$$

\noindent hence the surjectivity of $\Phi_{L}$ is implied by the
vanishing
of $H^1(S \times S,p_{1}^*L \otimes p_{2}^*L \otimes
{\cal I}_{\Delta}^2)
\cong H^1(Y,L_{1} \otimes L_{2}(-2E))$. On the other hand, since
$\omega_{Y} = {\cal O}_{Y}(E)$, we have $L_{1} \otimes L_{2}(-2E)
= \omega_{Y} \otimes
\bigotimes\limits_{j=1}^3 [A_{j1} \otimes A_{j2}(-E)]$ and the
required
vanishing follows by (ii) from the Kawamata-Viehweg vanishing
theorem.
\ \qed

As we will see below if we have three very ample line bundles as in
(i) of
Lemma (2.1) then (ii) is automatically satisfied. If at least one
line
bundle is very ample it is still possible that (ii) holds. Some
sufficient
conditions that are enough for our purposes are stated in the
following
lemma.

Let $A$ be a line bundle on a K3 surface $S$ with $A^2 \geq 2$,
having no
base points. Recall that we have the following three cases (see
[SD],
[Ma]) : if  $A^2 \geq 4, A \not= 2B$ with $B^2 = 2$ and the
associated
morphism $\phi_{A}$ is birational, then in fact it is an embedding
or an
isomorphism off some irreducible curves $Z$ such that $Z^2=-2, Z
\cdot A =
0$; if $A^2 = 2$ then $\phi_{A}$ is a 2:1 morphism onto $\I^2$
and is
finite if there are no irreducible curves $Z$ such that  $Z^2=-2,
Z \cdot A = 0$; if $A = 2B$ with $B^2 = 2$, then $\phi_{A}$ is a
2:1
morphism onto the Veronese surface in
$\I^5$. In these three cases we have

{\bf Lemma (2.2).} {\sl Let $A_{1}, A_{2}, A_{3}$ be three base
point
free line bundles on $S$ with $A_{j}^2 \geq 2, j=1,2,3$ and such
that
$A_{1}$ is very ample and $A_{2}, A_{3}$ define either isomorphisms
off
one (possibly empty) curve $Z_{2}, Z_{3}$ respectively or 2:1
finite
morphisms onto $\I^2$ or onto the Veronese surface in $\I^5$.
Suppose that
either

(i) $A_{2}$ and $A_{3}$ are very ample or

(ii) $(A_{1}+A_{2}+A_{3}) \cdot Z_{j} \geq 3$ whenever there is a
$Z_{j}$
and

(iii) $(A_{1}+A_{2}+A_{3}) \cdot A_{j} \geq 9$ whenever $A_{j}^2 =
2$ or
$(A_{1}+A_{2}+A_{3}) \cdot B_{j} \geq 9$ whenever $A_{j} = 2B_{j},
B_{j}^2 = 2$.

Then $\bigotimes\limits_{j=1}^3 [A_{j1} \otimes A_{j2} (-E)]$ is
big and
nef on $Y$.}

\noindent {\it Proof:} If $A_{j}$ is very ample the linear system
$|A_{j1} \otimes A_{j2} (-E)|$ on $Y$ has a sublinear system
defining the
morphism $Y \to $ {\ms G} $(1,\I H^0(A_{j})^*)$ associating to
$(x,y) \in Y$ the linear span of $\phi_{A_{j}}(x)$ and
$\phi_{A_{j}}(y)$
(note that this still makes sense if  $(x,y) \in E$ since we can
think of
$(x,y)$ as a pair with $x \in S$, $y \in \I T_{S_{|x}}$).
Therefore $A_{j1} \otimes A_{j2} (-E)$ is nef and also big since
the image
of $S$ in $\I H^0(A_{j})^*$ is not ruled. If $A_{2}$ and/or
$A_{3}$ is not
very ample certainly the line bundle $\bigotimes\limits_{j=1}^3
[A_{j1}
\otimes A_{j2} (-E)]$ is big since $A_{11} \otimes A_{12} (-E)$
is already
big; moreover it can fail to be nef only on a curve contained in
the
indeterminacy locus of the maps
$Y \to \hbox{\ms G}(1,\I H^0(A_{j})^*), j = 2,3$, that is a curve
$Z \subset \{(x,y) \in Y : \phi_{A_{j}}(x)=\phi_{A_{j}}(y) \}$. If
$A_{j}$
is an isomorphism off one curve $Z_{j}$ then $Z_{j} \cong \I^1$,
$Z \subset
Z_{j} \times Z_{j}$ (which we have identified with its blow-up in
$Y$), and
we  have
$$\bigotimes\limits_{j=1}^3 [A_{j1} \otimes A_{j2} (-E)]_{|Z_{j}
\times
Z_{j}}={\cal O}_{\I^1 \times \I^1}((A_{1}+A_{2}+A_{3}) \cdot
Z_{j},(A_{1}+A_{2}+A_{3}) \cdot Z_{j})(-3\Delta_{\I^1 \times
\I^1})=$$
$$= {\cal O}_{\I^1 \times \I^1}((A_{1}+A_{2}+A_{3}) \cdot
Z_{j}-3,(A_{1}+A_{2}+A_{3}) \cdot Z_{j}-3)$$

\noindent is nef on any curve $Z \subset Z_{j} \times Z_{j} \cong
\I^1
\times \I^1$ by (ii). Suppose now that $A_{j}$ gives a 2:1 finite
morphism
onto $\I^2$. If $Z \not\subset E$, letting $i : S \times S \to
S \times S$ be the  involution $i(x,y) = (y,x)$ we can assume that
the
image
$\overline Z$ of $Z$ in $S \times S$ is such  that $i(\overline Z)
=
\overline
Z$, because  $\bigotimes\limits_{j=1}^3 [A_{j1} \otimes A_{j2}(-E)]$
is
invariant under $i$. Then the first (or second) projection $Z'$ of
$\overline
Z$ in $S$ is $\phi_{A_{j}}^*(Z_{1})$ for some curve $Z_{1} \subset
\I^2$.
If
$N$ is a  line in $\I^2$ and $Z_{1} \sim mN$ we have
$$\bigotimes\limits_{j=1}^3 [A_{j1} \otimes A_{j2} (-E)] \cdot Z =
2(A_{1}+A_{2}+A_{3}) \cdot Z' -  3 E \cdot Z =$$ $$= 2
(A_{1}+A_{2}+A_{3})
\cdot m \phi_{A_{j}}^*(N) - 3 E \cdot Z = 2m (A_{1}+A_{2}+A_{3})
\cdot A_{j} -
3 E \cdot Z$$

\noindent and we will be done by (iii) if we show that $E \cdot Z
= 6m$.
To this end let $B$ be the ramification divisor of $\phi_{A_{j}}$;
then
$B$ is a smooth plane sextic and $E \cdot Z = mB \cdot N = 6m$
provided
that the intersection of $E$ and $Z$ is transverse. On the other
hand
there is a one to one correspondence between $E \cap Z$ and
$\Delta \cap
\overline Z$, and transversality can then be checked on $S
\times S$. For
our purpose it is of course enough to show that $\Delta$ and the
pull-back
of $N$ on $S
\times S$ intersect transversally. The latter being a local
computation,
we can assume that locally the double plane $S$ is given by
$z^2=x$ and
$N$ is $y=0$; then on $S \times S$ with coordinates $y,z,y',z'$
the
pull-back of $N$ is defined by the equations $z'=-z, y=y'=0$,
hence it is
transversal to the diagonal $\Delta$ of
$S \times S$.
If $Z \subset E$, since
$Z \subset \{(x,y) \in Y : \phi_{A_{j}}(x) = \phi_{A_{j}}(y) \}$,
then
it must be the strict transform of the ramification divisor
$\overline B$
on
$S$ of $\phi_{A_{j}}$, hence we have that
$Z \cdot E = - c_{1}({\cal O}_{\I T_{S}}(1)) \cdot Z = - deg \
T_{\overline
B} = 18$ and therefore
$\bigotimes\limits_{j=1}^3 [A_{j1} \otimes A_{j2} (-E)]
\cdot Z = 6 (A_{1}+A_{2}+A_{3}) \cdot A_{j} - 3 E \cdot Z \geq 0$
by (iii).
If
$A_{j} = 2B_{j}$ with $B_{j}^2 = 2$ then $\phi_{A_{j}} = v_{2}
\circ
\phi_{B_{j}}$, where $v_{2}$ is the Veronese map, and we have
$Z \subset
\{(x,y) \in Y : \phi_{A_{j}}(x) = \phi_{A_{j}}(y) \} = \{(x,y)
\in Y :
\phi_{B_{j}}(x) = \phi_{B_{j}}(y) \}$ hence on $S$,
$Z'= \phi_{A_{j}}^*(Z_{1})
= \phi_{B_{j}}^*(Z_{2})$ where  $Z_{1}$ is a curve on the Veronese
surface
in
$\I^5$ and $Z_{2} = v_{2}^*(Z_{1}) \subset \I^2$. Now setting
$Z_{2} \sim mN$
a computation as above shows that $$\bigotimes\limits_{j=1}^3
[A_{j1}
\otimes
A_{j2} (-E)] \cdot Z =  2(A_{1}+A_{2}+A_{3}) \cdot Z' - 3 E
\cdot Z =
2m(A_{1}+A_{2}+A_{3}) \cdot B_{j} - 18m$$

\noindent and we are done by (iii) (similarly when $Z \subset E$).
\ \qed

We will now start the proof of Theorem (1.1). This will be done in
two
main steps: First we will compute the corank of
$\Phi_{{\cal O}_{S_{i,g}}(1)}$ for low values of $g$ or high values
of $i$
and then (in section 3, for $i=1, g
\geq 17$ or $i=2, g \geq 7$) we will construct K3 surfaces in
${\cal H}_{g}$ whose {\it i}-th Veronese embedding has the
hyperplane
bundle
decomposable as in Lemma (2.1), and hence
$\Phi_{{\cal O}_{S_{i,g}}(1)}$
surjective.

\noindent {\it Proof of Theorem (1.1):} First of all, by the
semicontinuity
of the corank of Gaussian maps, in order to prove that the general
surface
in an irreducible family has a surjective Gaussian map, we only
need to
exhibit a single surface which has a surjective Gaussian map.
For $i \geq 3$ and $g \geq 3$ on any smooth K3 surface representing
a point
in ${\cal H}_{i,g}$ there is an obvious decomposition of the
hyperplane
bundle as a tensor product of three very ample line bundles, hence
satisfying (ii) of Lemma (2.1) by Lemma (2.2), and therefore giving
the
required surjectivity by Lemma (2.1). For $g = 2$ the values of
$corank \
\Phi_{{\cal O}_{S_{i,2}}(1)}$ follow by a result of J. Duflot
[D,
Proposition 4.7]. Alternatively, for $g = 2, i \geq 5$ the
surjectivity can
also be proved using Lemmas (2.1) and (2.2) as follows.  We have
$iH = 3H +
H + (i-4)H$ and the conditions of Lemma (2.2) are satisfied with
$A_{1} =
3H, A_{2} = H, A_{3} = (i-4)H$: In fact $3H$ is very ample and
$(i-4)H$ is
very ample for $i \geq 7$ while it defines a 2:1 finite morphism
for $i =
5, 6$ and we have $(A_{1}+A_{2}+A_{3}) \cdot H = 2i \geq 10$. For
$i = 1, g
\geq 17$ or  $i = 2, g \geq 7$ there is a K3 surface $T$ whose
{\it i}-th
Veronese embedding $v_{i}(T)$ has surjective Gaussian map by
Proposition
(3.1) and Lemma (2.1).

Now for $L = {\cal O}_{S_{i,g}}(1)$ set $R(L,L^k) = Ker \{ H^0(L)
\otimes
H^0(L^k) \to H^0(L^{k+1}) \}$ and consider the Gaussian map
$$\Phi_{L,L^k} : R(L,L^k) \to H^0(S_{i,g},\Omega_{S_{i,g}}^1
\otimes
L^{k+1})$$
\noindent defined as usual by $\Phi_{L,L^k} (\sigma \otimes \tau) =
\sigma
d\tau - \tau d\sigma$. We recall that for $k = 1$ we have $Im \
\Phi_{L,L}
= Im \ \Phi_{L}$, where $\matrix{\Phi_{L} : \bigwedge^2
H^0(S_{i,g},L)  \to H^0(S_{i,g},\Omega_{S_{i,g}}^1 \otimes L^2)
\cr}$. With the
same notation as in Lemma (2.1) we  have
$$Coker \Phi_{L,L^k} \subseteq H^1(S \times S,p_{1}^*L \otimes
p_{2}^*L^k
\otimes {\cal I}_{\Delta}^2) \cong H^1(Y, \omega_{Y} \otimes L_{1}
\otimes
L_{2}^k(-3E)) = 0$$ \noindent by the Kawamata-Viehweg vanishing
theorem for
$i=1$ and $g \geq 17, i=2$ and $g \geq 7, i \geq 3$ and $g \geq 3,
g = 2$ and
$i \geq 5$. In fact as we have seen above in all these cases a
decomposition
as in (i) and (ii) of Lemma (2.1) holds for $L$ so that $L_{1}
\otimes
L_{2}(-3E) = \bigotimes\limits_{j=1}^3 [A_{j1} \otimes A_{j2} (-E)]$
is big
and nef on $Y$, and therefore so is  $L_{1} \otimes L_{2}^k(-3E)$.

It remains to prove the surjectivity of $\Phi_{L,L^k}$ for $k \geq 1,
i = 2, g
= 3, 4, 5, 6$. Since $S_{2,g} = v_{2}(S_{g})$, where $v_{2}$ is the
second
Veronese embedding, we clearly have $corank \ \Phi_{L,L^k} =
corank \ \Phi_{{\cal O}_{S_{g}}(2),{\cal O}_{S_{g}}(2k)}$. When $g
\leq 5$ the
K3 surface $S_{g}$ is a complete intersection and the corank of the
Gaussian
map $\Phi_{{\cal O}_{S_{g}}(2),{\cal O}_{S_{g}}(2k)}$ can be
computed by
results of Wahl and S. Kumar as follows. We have a diagram
$$\matrix{R({\cal O}_{\I^g}(2),{\cal O}_{\I^g}(2k))  &
\mapright{\Phi_{{\cal
O}_{\I^g}(2),{\cal O}_{\I^g}(2k)}} & H^0(\I^g,\Omega_{\I^g}^1(2k+2))
\cr    &&
\downarrow \phi_{g} \cr  \downarrow \pi_{g} & & \ \
H^0(S_{g},\Omega_{\I^g}^1(2k+2)_{|S_{g}}) \cr   &&\downarrow \psi_{g}
\cr
R({\cal O}_{S_{g}}(2),{\cal O}_{S_{g}}(2k)) &\mapright{\Phi_{{\cal
O}_{S_{g}}(2),{\cal O}_{S_{g}}(2k)}}  & H^0(S_{g},
\Omega_{S_{g}}^1(2k+2)) \cr}
\leqno (2.3)$$
\noindent and $\Phi_{{\cal O}_{\I^g}(2),{\cal O}_{\I^g}(2k)}$ is
 surjective by
[W1, Theorem 6.4] and [K, Theorem 2.5]. Since $S_{g} \subset
\I^{g}$ is a
complete intersection surface of type $(4)$ for $g = 3$, type
$(2,3)$ for  $g =
4$ and type $(2,2,2)$ for $g = 5$, we have that $N_{S_{g}
/\I^{g}}^*(2k+2) =
{\cal O}_{S_{3}}(2k-2), {\cal O}_{S_{4}}(2k) \oplus
{\cal O}_{S_{4}}(2k-1),
{\cal O}_{S_{5}}(2k)^{\oplus 3}$ respectively and its $H^1$
vanishes
(in fact $H^1({\cal O}_{S_{g}}(a)) = 0$ for every integer $a$);
therefore
$\psi_{g}$ is surjective since $Coker \psi_{g} \subseteq
H^1(N_{S_{g}
/\I^{g}}^*(2k+2)) = 0$. Now $Coker \phi_{g} =
H^1(\Omega_{\I^g}^1
\otimes {\cal I}_{S_{g} /\I^{g}}(2k+2))$ since
$H^1(\Omega_{\I^g}^1(2k+2)) = 0$
by Bott vanishing. For $g = 4, 5$ and the Koszul resolution of
the ideal
sheaf ${\cal I}_{S_{g} /\I^{g}}$ we have
$$0 \to {\Omega_{\I^4}^1}(2k-3) \to
{\Omega_{\I^4}^1}(2k) \oplus {\Omega_{\I^4}^1}(2k-1) \to
{\Omega_{\I^4}^1}
\otimes {\cal I}_{S_{4} /\I^{4}}(2k+2) \to 0$$
\noindent and
$$0 \to {\Omega_{\I^5}^1}(2k-4) \to
{\Omega_{\I^5}^1(2k-2)}^{\oplus 3}
\to {\Omega_{\I^5}^1}(2k)^{\oplus 3} \to {\Omega_{\I^5}^1}
\otimes {\cal
I}_{S_{5} /\I^{5}}(2k+2) \to 0$$
\noindent and by Bott vanishing we see that
$H^1(\Omega_{\I^g}^1 \otimes {\cal I}_{S_{g} /\I^{g}}(2k+2)) = 0$,
hence
$\phi_{g}$ is surjective and so is $\Phi_{{\cal O}_{S_{g}}(2),{\cal
O}_{S_{g}}(2k)}$, by  diagram (2.3).

\noindent For $g = 3$ we have that
$corank \ \phi_{3} = h^1(\Omega_{\I^3}^1(2k-2)) = 0$ unless $k=1$.
In the
latter case we get $corank \ \phi_{3} = 1$ and  $dim Ker
\psi_{3} = h^0({\cal O}_{S_{3}}) = 1$; hence to see the  surjectivity
of
$\Phi_{{\cal O}_{S_{3}}(2)}$ it is enough to show that  $Ker \psi_{3}
\not\subseteq Im \phi_{3}$ because then $\psi_{3} \circ \phi_{3}$ is
surjective. Now consider the diagram
$$\matrix{0 & \to
H^0({\Omega_{\I^3}^1}(4)) & \to &  H^0({\cal O}_{\I^3}(1)) \otimes
H^0({\cal
O}_{\I^3}(3)) & \to & H^0({\cal O}_{\I^3}(4)) & \to & 0 \cr & \
\downarrow{\phi_{3}} && \ \ \ \downarrow{\cong} &&\downarrow{\alpha}
\cr 0 &
\to H^0({\Omega_{\I^3}^1}_{|S_{3}}(4)) & \to &
H^0({\cal O}_{S_{3}}(1)) \otimes
H^0({\cal O}_{S_{3}}(3)) & \to &  H^0({\cal O}_{S_{3}}(4)) & \to &
0 \cr}$$

\noindent and let $\delta : H^0({\cal I}_{S_{3} /\I^{3}}(4)) = Ker
\alpha \to Coker \phi_{3}$ be the isomorphism induced by the snake
lemma.
Since $\delta$ is surjective and $Ker \psi_{3} = H^0(N_{S_{3}
/\I^{3}}^*(4)) = H^0({\cal I}_{S_{3} /\I^{3}}(4))$, to show that
$Ker
\psi_{3}  \not\subseteq Im \phi_{3}$ is enough to see that $4
\delta$ is
the
composition $H^0(N_{S_{3} /\I^{3}}^*(4))
\to H^0({\Omega_{\I^3}^1}_{|S_{3}}(4)) \to Coker \phi_{3}$. To this
end
observe that if $F = 0$ is the equation of $S_{3}$, we have $F =
{1 \over 4} \sum\limits_{i=0}^3 {\partial F \over \partial x_{i}}
x_{i}$
hence by the snake lemma applied to the above diagram, $4
\delta(F) =
\sum\limits_{i=0}^3 x_{i} \otimes {\partial F \over \partial x_{i}}
+ Im
\phi_{3}$ (where we see $\sum\limits_{i=0}^3 x_{i} \otimes
{\partial F
\over \partial x_{i}}$ as an element of the kernel of the
multiplication
map $\mu : H^0({\cal O}_{S_{3}}(1)) \otimes
H^0({\cal O}_{S_{3}}(3)) \to
H^0({\cal O}_{S_{3}}(4))$, hence as an element of
$H^0({\Omega_{\I^3}^1}_{|S_{3}}(4))$). On the other hand the
map
$H^0(N_{S_{3} /\I^{3}}^*(4)) \to
H^0({\Omega_{\I^3}^1}_{|S_{3}}(4))$
takes $F$ to $dF = \sum\limits_{i=0}^3 {\partial F \over
\partial x_{i}}
dx_{i}$ which viewed as an element of the kernel of $\mu$ is in fact
$\sum\limits_{i=0}^3 x_{i} \otimes {\partial F \over \partial
x_{i}}$.

Finally to see the surjectivity of
$\Phi_{{\cal O}_{S_{6}}(2),{\cal
O}_{S_{6}}(2k)}$ we will use the fact that $S_{6}$ is a complete
intersection
in a Grassmannian. We have $\hbox{\ms G} = \hbox{\ms G}(1,4)
\subset \I^{9}$ in
the Pl\"ucker embedding and $S_{6} = \hbox{\ms G} \cap H_{1}
\cap H_{2} \cap
H_{3} \cap Q$ where $H_{i}$ is a hyperplane and $Q$ is a quadric
hypersurface.
In the diagram
$$\matrix{R({\cal O}_{\hbox{\ms G}}(2),{\cal O}_{\hbox{\ms G}}(2k))
 &
\mapright{\Phi_{{\cal O}_{\hbox{\msp G}}(2),
{\cal O}_{\hbox{\msp G}}(2k)}} &
H^0(\hbox{\ms G},\Omega_{\hbox{\ms G}}^1(2k+2)) \cr
&& \downarrow \phi_{6}
\cr  \downarrow \pi_{6} & & \ \ H^0(S_{6},\Omega_{\hbox{\ms
G}}^1(2k+2)_{|S_{6}}) \cr   &&\downarrow \psi_{6} \cr R({\cal
O}_{S_{6}}(2),{\cal O}_{S_{6}}(2k)) &\mapright{\Phi_{{\cal
O}_{S_{6}}(2),{\cal
O}_{S_{6}}(2k)}}  & H^0(S_{6},  \Omega_{S_{6}}^1(2k+2)) \cr} $$
\noindent we have that $\Phi_{{\cal O}_{\hbox{\msp G}}(2),{\cal
O}_{\hbox{\msp
G}}(2k)}$ is surjective by [W1, Theorem 6.4] and [K, Theorem 2.5],
$\psi_{6}$ is surjective because $Coker \psi_{6} \subseteq
H^1(N_{S_{6}
/\hbox{\msq G}}^*(2k+2)) = H^1({\cal O}_{S_{6}}(2k+1)^{\oplus 3}
\oplus {\cal
O}_{S_{6}}(2k)) = 0$. Therefore the surjectivity of
$\Phi_{{\cal O}_{S_{6}}(2),{\cal O}_{S_{6}}(2k)}$ will follow by
the above
diagram as soon as we show that $\phi_{6}$ is surjective. To this
end note
that $Coker \phi_{6} \subseteq H^1(\Omega_{\hbox{\ms G}}^1 \otimes
{\cal
I}_{S_{6} /\hbox{\msq G}}(2k+2)) = 0$: by the Koszul resolution of
the ideal
sheaf ${\cal I}_{S_{6} /\hbox{\msq G}}$ we have an exact sequence
$$0 \to {\Omega_{\hbox{\ms G}}^1}(2k-3) \to {\Omega_{\hbox{\ms
G}}^1}(2k-2)^{\oplus 3} \oplus {\Omega_{\hbox{\ms G}}^1}(2k-1) \to
{\Omega_{\hbox{\ms G}}^1}(2k)^{\oplus 3} \oplus {\Omega_{\hbox{\ms
G}}^1}(2k-1)^{\oplus 3} \to$$
$$ \to {\Omega_{\hbox{\ms G}}^1}(2k+1)^{\oplus 3} \oplus
{\Omega_{\hbox{\ms G}}^1}(2k) \to {\Omega_{\hbox{\ms G}}^1} \otimes
{\cal
I}_{S_{6} / \hbox{\msq G}}(2k+2) \to 0$$
\noindent and the vanishing follows by the

\noindent {\it Claim} (2.4). {\sl $H^p(\Omega_{\hbox{\ms G}}^1(q))
= 0 \
\hbox{for} \  p \geq 2, q \geq -1 \ \hbox{and for} \ p=1, q\geq 1$.}

\noindent {\it Proof of Claim (2.4):} This is just Bott vanishing
for the
Grassmannian. Alternatively from the normal bundle sequence
$$0 \to N_{\hbox{\ms G}}^*(q) \to
\Omega_{\I^{9}}^1(q)_{|\hbox{\ms G}}
\to \Omega_{\hbox{\ms G}}^1(q) \to 0$$
\noindent we see that the Claim is implied by
$$H^p(\Omega_{\I^{9}}^1(q)_{|\hbox{\ms G}}) = 0 \ \hbox{for} \
p \geq 2, q \geq -3 \ \hbox{and for} \ p=1, q \geq 1 \leqno (2.5)$$
\noindent and
$$H^{p+1}(N_{\hbox{\ms G}}^*(q)) = 0 \ \hbox{for} \ p \geq 0, q \geq
-1. \leqno
(2.6)$$
\noindent To see (2.5) consider the Euler sequence
$$0 \to \Omega_{\I^{9}}^1(q)_{|\hbox{\ms G}} \to H^0({\cal
O}_{\hbox{\ms G}}(1)) \otimes {\cal O}_{\hbox{\ms G}}(q-1) \to
{\cal O}_{\hbox{\ms G}}(q) \to 0.$$
\noindent We have $H^p({\cal O}_{\hbox{\ms G}}(q - 1)) =
H^p(\omega_{\hbox{\ms G}}(4 + q)) = 0$ for $p \geq 1, q \geq -3$ by
Kodaira
vanishing and $H^{p-1}({\cal O}_{\hbox{\ms G}}(q)) =
H^{p-1}(\omega_{\hbox{\ms G}}(5 + q)) = 0$ for $p \geq 2, q \geq -4$
again by
Kodaira vanishing, therefore (2.5) holds for $p \geq 2, q \geq -3$
and also for
$p = 1, q \geq 1$ since by what we just proved we have
$$H^0({\cal O}_{\hbox{\ms G}}(1)) \otimes
H^0({\cal O}_{\hbox{\ms G}}(q-1)) \to
H^0({\cal O}_{\hbox{\ms G}}(q)) \to
H^1(\Omega_{\I^{9}}^1(q)_{|\hbox{\ms G}}) \to  0$$
\noindent and the above multiplication map is surjective.

\noindent Now to prove (2.6) we use Griffiths vanishing theorem
([G]).
In the Pl\"ucker embedding the ideal of the Grassmannian {\ms G} is
generated by quadrics, hence $N_{\hbox{\ms G}}^*(2)$ is globally
generated
and $det (N_{\hbox{\ms G}}^*(2)) = det (N_{\hbox{\ms G}}^*)(6) =
\omega_{\hbox{\ms G}}^{-1}(-4) =  {\cal O}_{\hbox{\ms G}}(1)$.
Therefore
setting $E = N_{\hbox{\ms G}}^*(2)$ we can write
$$N_{\hbox{\ms G}}^*(q) = E \otimes det E \otimes {\cal
O}_{\hbox{\ms
G}}(q+2) \otimes \omega_{\hbox{\ms G}}$$

\noindent with $p \geq 0, q + 2 > 0$ and (2.6) follows by Griffiths
vanishing
theorem. \ \qed \ (for Claim (2.4))

\noindent This then concludes the proof of Theorem (1.1). \ \qed

{\bf (2.7)} {\it Remark.} Note that in the proof of Theorem (1.1)
for $i=2$
and $3 \leq g \leq 4$ or for $i \geq 3$ and $g \geq 3$ {\it we did
not use} the
fact that $S_{i,g}$ represents a general point of ${\cal H}_{i,g}$.
Hence for
such values of $i$ and $g$ the corank of
$\Phi_{{\cal O}_{S_{i,g}}(1)}$ and
of $\Phi_{{\cal O}_{S_{i,g}}(1),{\cal O}_{S_{i,g}}(k)}$ is as in
Theorem (1.1) for {\it any smooth} $S_{i,g}$. Similarly in the case
$g=2$ we
only use the fact that the ramification locus is smooth (as in [D,
Proposition
4.7]) hence the computation of the corank of
$\Phi_{{\cal O}_{S_{i,2}}(1)}$
holds just with this hypothesis.

\vskip .3cm
\noindent {\bf 3. CONSTRUCTION OF K3 SURFACES WITH SURJECTIVE
GAUSSIAN
MAP}
\vskip .3cm
In this section we will construct with the aid of the surjectivity
of the
period mapping, K3 surfaces with very ample line bundles whose
multiples
decompose as in (ii) of Lemma (2.1). We have

{\bf Proposition (3.1).} {\sl For every $i,g$ such that $i=1, g
\geq 17$ or
$i=2, g \geq 7$ there exists a K3 surface $T$ representing a point
in
${\cal H}_{g}$ such that $v_{i}(T)$ has the  hyperplane bundle
decomposable as in Lemma (2.1).}

The proof of Proposition (3.1) will be in two parts. We will first
construct the needed K3 surfaces using the surjectivity of the
period map
and then we will show how to decompose $iH, i=1,2$.

Let us recall that a cohomology class $h \in H^2(X,\hbox{\ms R})$ on
a
compact K\"ahler manifold $X$ is a {\it K\"ahler class} if it can
be
represented by a K\"ahler form, that is by a (1,1)-form associated
to a
K\"ahler metric. In particular if $h$ is a K\"ahler class and $y$
is the
class of any closed curve on
$X$ we have $h \cdot y > 0$ [BPV, Lemma I.13.1].

{\bf Proposition (3.2).} {\sl There exist smooth K3 surfaces
$T_{jkh}$
with Picard lattices $\Gamma_{jkh}$ as follows:

\noindent (i) $\Gamma_{jkh} =$ {\Z}$D \oplus ${\Z}$L$ with
intersection
matrix $\left(\matrix{D^2 & D \cdot L \cr L \cdot D & L^2
\cr}\right) =
\left(\matrix{2h & k \cr k & 2j \cr}\right)$

\noindent for $j=1,2, k \geq j+4, h=2$, or $j=-1, k=1,2, h
\geq 5-2k$ or
$j=1, k=5, h=3$;

\noindent (ii) $\Gamma_{jkh} =$ {\Z}$D \oplus ${\Z}$L \oplus
${\Z}$R$
with intersection matrix
$$\left(\matrix{D^2 & D \cdot L & D \cdot R \cr
L \cdot D & L^2 & L \cdot R \cr R \cdot D & R \cdot L & R^2
\cr}\right)
= \left(\matrix{2h & k & 2 \cr k & 4j-2 & j \cr 2 & j & -2
\cr}\right)$$
\noindent for $j=0, k=1,2, h \geq 5-2k$ and $j=1, k=4, h=1$.

\noindent Moreover all the $D$'s are K\"ahler (hence ample)
classes on
$T_{jkh}$.}

\noindent {\it Proof:} Let $\Lambda = \hbox{\ms H}^3 \oplus
E_{8}(-1)^2$
be the K3 lattice where {\ms H} denotes the hyperbolic plane
and $E_{8}$
the root lattice, $\Lambda_{\hbox{\ms C}} = \Lambda \otimes
{\hbox{\ms
C}}$, $\Lambda_{\hbox{\ms R}} = \Lambda \otimes {\hbox{\ms R}}$
with the
quadratic form extended {\ms C} or {\ms R}-bilinearly and let
$$(K\Omega)^0 = \{(k,[\omega]) \in \Lambda_{\hbox{\ms R}} \times
\I\Lambda_{\hbox{\ms C}} : \omega \cdot \omega = 0, \omega \cdot
\overline{\omega} > 0, k \cdot k > 0, k \cdot \omega = 0 \
\hbox{and} \
k \cdot d \not= 0$$

\ \ \ \ \ \ \ \ \ $\forall d \in \Lambda$ with $d^2 = -2, \omega
\cdot
d = 0 \}.$

We will show that there exists $\omega \in \Lambda_{\hbox{\ms C}}$
such
that $(D,[\omega]) \in (K\Omega)^0$. To this end let us record the
following

\noindent {\it Claim} (3.3). {\sl The lattices $\Gamma_{jkh}$ are
even,
nondegenerate of signature $(1, rk \Gamma_{jkh} - 1)$ and there is
no
class $d \in \Gamma_{jkh}$ such that $d^2 = -2, D \cdot d = 0$.}

\noindent {\it Proof of the Claim:} Clearly all the $\Gamma_{jkh}$
are
even and nondegenerate since $disc \Gamma_{jkh} =$

\noindent $det$ (intersection matrix) $\not= 0$
(see Table (3.4)). The assertion about the signatures follows easily
by
looking at the signs of the principal minors of the intersection
matrix. Now suppose there is  $d \in \Gamma_{jkh}$ such that $d^2 =
-2,  D
\cdot d = 0$. In case (i), since $disc \Gamma_{jkh}$ must
\eject
\noindent divide the discriminant of any rank two sublattice, we get
that
$4hj-k^2$ divides $disc(D,d) = -4h$ and this does not hold for the
given
values of $j,k,h$. In case (ii) set $d = \alpha D + \beta L +
\gamma R$; then
$d \cdot D = 0$ gives $\gamma = - h \alpha - {k \over 2} \beta$ and
from $d^2
= -2$ we get $$2h \alpha^2 + (4j-2)\beta^2-2(- h \alpha - {k \over 2}
\beta)^2
+ 2k \alpha \beta + 4 \alpha (- h \alpha - {k \over 2} \beta) + 2j
\beta (- h
\alpha - {k \over 2} \beta) = -2$$

\noindent hence
$$2h(h+1)\alpha^2+2h\beta(k+j)\alpha + (jk + {k^2 \over 2} + 2 -
4j)\beta^2 - 2 = 0$$

\noindent whose discriminant in $\alpha$ is
$$\Delta_{\alpha} = h[h(k+j)^2-2(h+1)(jk + {k^2 \over 2} + 2 -
4j)]\beta^2 + 4h(h+1),$$

\noindent which is easily seen to be negative for $\beta \not= 0$,
while
for $\beta = 0$ it is not a square since $4h(h+1)$ is not a square
for $h
\geq 1$. \ \qed \ (for Claim (3.3))

To finish the proof of Proposition (3.2) let $\Gamma =
\Gamma_{jkh}$.
Since $rk \Gamma \leq 3$, by [BPV, Theorem I.2.9], $\Gamma$ has a
primitive embedding into the K3 lattice $\Lambda$ and by Claim (3.3)
the
signature of $\Gamma^{\perp}$ is $(2, 20 - rk \Gamma)$. Hence there
is a
positive definite two dimensional space $V \subset
\Gamma_{\hbox{\ms
R}}^{\perp}$ such that $V^{\perp} \cap \Lambda = \Gamma$. Indeed
let $u$
and $v$ be two orthogonal vectors in $\Gamma_{\hbox{\ms R}}^{\perp}$
spanning such a positive definite two dimensional space; multiplying
by a
real factor we can assume $u^2=v^2$, hence if we set $\omega_{1} =
u +iv
\in \Lambda_{\hbox{\ms C}}$ we have $\omega_{1}^2 = 0$, $\omega_{1}
\cdot
\overline{\omega_{1}} > 0$ and $\Gamma \subseteq (\hbox{\ms
C}\omega_{1}
\oplus \hbox{\ms C}\overline{\omega_{1}})^{\perp}\cap \Lambda$.
By [EEK,
proof of Theorem 5.4], $\omega_{1}$ can be perturbed, to a class
$\omega
\in \Lambda_{\hbox{\ms C}}$, in such a way as to preserve the
first two
relations and achieve equality in the third. Then we can take $V =
\hbox{\ms R}Re \omega \oplus \hbox{\ms R}Im \omega$. Now let $d
\in
\Lambda$ be such that $d^2 = -2, \omega \cdot d = 0$. Then we get
that
$d \in V^{\perp}$. In fact if $d = d_{1} + d_{2}, \ d_{1} \in V,
\ d_{2} \in V^{\perp}$, we have $d_{1} = a Re \omega + b Im
\omega$ and
$0 = \omega \cdot d = \omega \cdot d_{1} = (Re \omega + i Im
\omega)
\cdot (a Re \omega + b Im \omega) = a (Re \omega)^2 + ib (Im
\omega)^2$
hence $a = b = 0$, that is $d_{1} = 0$. Therefore $d \in V^{\perp}
\cap
\Lambda = \Gamma$ and by Claim (3.3) we have $D \cdot d \not= 0$.
Since
$D^2 > 0$ we have that $(D,[\omega]) \in (K\Omega)^0$ and hence,
by the
surjectivity of the refined period mapping [BPV, Theorem VIII.1.4]
there exists a marked K3 surface $T_{jkh}$ with period point
$[\omega]$,
Picard lattice $\Gamma = \Gamma_{jkh}$ and K\"ahler class $D$.
The latter
implies that $D$ is an ample class. \ \qed

\noindent {\it Proof of Proposition (3.1):} On the K3 surfaces
$T_{jkh}$
with Picard lattices $\Gamma_{jkh}$ given in Proposition (3.2) let
us
define divisors $H, A_{1}, A_{2}, A_{3}$ according to the following
table
(in which we record also the genus $g(H) = {1 \over 2}H^2 + 1$ and
the
discriminants):

\centerline{\bf Table (3.4)}

\vbox{\offinterlineskip
\hrule
\halign{&\vrule#&
\strut\hfil#\tabskip=.28em\cr
height4pt&\omit&&\omit&&\omit&&\omit&&\omit&&\omit&&\omit&&\omit
&&\omit&&\omit&&\omit&\cr
&\ \ i\hfil&&\ j\hfil&&\ k\hfil&&\ h\hfil&&\ \ \ rank\hfil&&\ \
$disc\Gamma_{jkh}$\hfil &&\ \ $H$\hfil&&\ \ $g(H)$\hfil&&\ \
$A_{1}$\hfil&&\ $A_{2}$\hfil&&$A_{3}$\hfil&\cr
height4pt&\omit&&\omit&&\omit&&\omit&&\omit&&\omit&&\omit&&\omit
&&\omit&&\omit&&\omit&\cr
\noalign{\hrule}
height4pt&\omit&&\omit&&\omit&&\omit&&\omit&&\omit&&\omit&&\omit
&&\omit&&\omit&&\omit&\cr
&1&&1,2&&$\geq j+4$\hfil&&2&&2\hfil&&$8j-k^2$\hfil&&$2D+L$\hfil
&&\ \
$2k+j+9$\hfil&&$D$\hfil&&$D$\hfil&&$L$\hfil&\cr
height4pt&\omit&&\omit&&\omit&&\omit&&\omit&&\omit&&\omit&&\omit
&&\omit&&\omit&&\omit&\cr
&1&&1&&5,6,7\hfil&&2&&2\hfil&&$8-k^2$\hfil&&$D+2L$\hfil &&\ \
$2k+7$\hfil&&$D$\hfil&&$L$\hfil&&$L$\hfil&\cr
height4pt&\omit&&\omit&&\omit&&\omit&&\omit&&\omit&&\omit&&\omit
&&\omit&&\omit&&\omit&\cr
&1&&1&&5\hfil&&3&&2\hfil&&$-13$\hfil&&$D+2L$\hfil &&\ \
18\hfil&&$D$\hfil&&$L$\hfil&&$L$\hfil&\cr
height4pt&\omit&&\omit&&\omit&&\omit&&\omit&&\omit&&\omit&&\omit
&&\omit&&\omit&&\omit&\cr
&2&&$-1$&&1\hfil&&\ $\geq 3$&&2\hfil&&$-4h-1$\hfil&&$2D-L$\hfil
&&\ \ $4h-2$\hfil&&$H$\hfil&&$D$\hfil&&$D-L$\hfil&\cr
height4pt&\omit&&\omit&&\omit&&\omit&&\omit&&\omit&&\omit&&\omit
&&\omit&&\omit&&\omit&\cr
&2&&$-1$&&2\hfil&&$\geq 1$&&2\hfil&&$-4h-4$\hfil&&$2D+L$\hfil
&&\ \ $4h+4$\hfil&&$H$\hfil&&$D$\hfil&&$D+L$\hfil&\cr
height4pt&\omit&&\omit&&\omit&&\omit&&\omit&&\omit&&\omit&&\omit
&&\omit&&\omit&&\omit&\cr
&2&&0&&1\hfil&&$\geq 3$&&3\hfil&&$8h+10$\hfil&&\ $2D-L+R$\hfil
&&\ \ $4h+1$\hfil&&$H$\hfil&&$D$\hfil&&\ $D-L+R$\hfil&\cr
height4pt&\omit&&\omit&&\omit&&\omit&&\omit&&\omit&&\omit&&\omit
&&\omit&&\omit&&\omit&\cr
&2&&0&&2\hfil&&$\geq 1$&&3\hfil&&$8h+16$\hfil&&\ $2D+L+R$\hfil
&&\ \ $4h+7$\hfil&&$H$\hfil&&\ $D+L$\hfil&&$D+R$\hfil&\cr
height4pt&\omit&&\omit&&\omit&&\omit&&\omit&&\omit&&\omit&&\omit
&&\omit&&\omit&&\omit&\cr
&2&&1&&5\hfil&&2&&2\hfil&&$-17$\hfil&&$D+L$\hfil &&\
\ 9\hfil&&$H$\hfil&&$D$\hfil&&$L$\hfil&\cr
height4pt&\omit&&\omit&&\omit&&\omit&&\omit&&\omit&&\omit&&\omit
&&\omit&&\omit&&\omit&\cr
&2&&1&&4\hfil&&1&&3\hfil&&30\hfil&&$D+L$\hfil &&\ \ 7\hfil&&$H$\hfil
&&$D$\hfil&&$L$\hfil&\cr
height4pt&\omit&&\omit&&\omit&&\omit&&\omit&&\omit&&\omit&&\omit
&&\omit&&\omit&&\omit&\cr
} \hrule}

Note that the first three rows of the above table are relative to
the
case $i=1, g \geq 17$, while the remaining six to the case $i=2, g
\geq 7$.

We will show that $H$ is very ample and embeds $T_{jkh}$ as a K3
surface
$T \subset \I^g$ with the required properties. To see this first
observe
that $g(H)$ runs through all the integers $g \geq 17$ for $i=1$ and
$g \geq 7$ for $i=2$ as Table (3.4) shows. Also notice that in all
cases
$H$ is a linear combination of generators of the Picard group in
which
appears at least one generator with coefficient $\pm1$, hence
$H$ is
indivisible and, by the transcendental theory of K3 surfaces, $T$
represents a point in
${\cal H}_{g}$. To finish the proof we will see, with a case by
case
analysis, that $H$ is very ample and that the decomposition of
$iH, i=1,2,$
given in Table (3.4) satisfies (ii) of Lemma (2.1).

Before we start we record, for the reader's convenience, the
well-known
theorems that we will use.

{\bf Lemma (3.5).} {\sl Let $S$ be a smooth K3 surface and
$\Delta$ an
indivisible divisor on $S$ with $\Delta^2 \geq 2$. We have

(i) $\Delta$ is very ample if $\Delta$ is nef, $\Delta^2 \geq 4$
and
there is no irreducible curve $F$ on $S$ such that either $F^2 = 0,
F \cdot \Delta =1,2$ or $F^2 = -2, F \cdot \Delta = 0$;

(ii) If $\Delta$ is nef it has no base points unless
there exist irreducible curves $F$, $G$ and an integer $a \geq 2$
such
that $\Delta \sim aF + G$, $F^2 = 0, G^2 = -2, F \cdot G =1$;

(iii) If $\Delta^2 = 2$ and $\Delta$ has no base points then it
defines
a 2:1 morphism $S \to \I^2$. Moreover the morphism is finite if
there is
no irreducible curve $F$ on $S$ such that $F^2 = -2, F \cdot
\Delta = 0$;

(iv) If $\Delta^2 \geq 4$, $\Delta$ has no base points and there is
no
irreducible curve $F$ on $S$ such that $F^2 = 0, F \cdot \Delta =
2$
then $\Delta$ defines a birational morphism. When the latter holds
the
only curves contracted by $\Delta$ are the curves $F$ such that
$F^2 = -2, F \cdot \Delta = 0$.}

\noindent {\it Proof of Lemma (3.5):} By Mori's version of a theorem
of
Saint-Donat [Mo, Theorem 5] it follows that under the hypotheses of
case
(i), if $\Delta$ is  not very ample then there is an irreducible
curve $F$
on $S$ such that $\Delta \sim 2F$, hence $\Delta$ is divisible.
Part (ii)
follows again by [Mo, Theorem 5]; (iii) is a result of Mayer [Ma,
Proposition 2] and  Saint-Donat ([SD]). By a theorem of Saint-Donat
([SD])
under the hypotheses of case (iv), if $\Delta$ does not define a
birational morphism then there is an irreducible curve $F$ such
that
$\Delta \sim 2F$, hence $\Delta$ is divisible. The rest follows
again by
[SD]. \ \qed \ (for Lemma (3.5))

\noindent {\it Claim} (3.6). {\sl For $j= -1, k=1, h \geq 3$ we have
that
$D$ is very ample and so is $D-L$ unless $h=3$ and in that case
$D-L$
defines a 2:1 finite morphism onto $\I^2$.}

By Claim (3.6) we see that, for $j=-1, k=1, h \geq 3$, $H = A_{1} =
D + D
- L$ is very ample and (ii) of Lemma (2.1) holds by Lemma (2.2)
since for
$h \geq 4$ $A_{2}$ and $A_{3}$ are very ample and for $h=3$ we have
$(A_{1}+A_{2}+A_{3}) \cdot A_{3} = (4D - 2L) \cdot (D - L) = 14.$

\noindent {\it Proof of Claim (3.6):} By Proposition (3.2) $D$ is a
K\"ahler (hence ample) indivisible class; if there is an irreducible
curve
$F$ such that $F^2 = 0, F \cdot D =1,2$ ($F \cdot D = 0$ is excluded
since
$D$ is an ample class) then $disc \Gamma_{-1,1,h} = -4h-1$ divides
$disc(F,D) =-1,-4$. Hence $D$ is very ample by (i) of Lemma (3.5). If
$h
\geq 4$ we have $(D - L)^2 = 2h - 4 \geq 4$ and $D - L$ is indivisible
and
nef since $(D - L) \cdot L = 3$  ($L$ is effective and irreducible
because
$L^2 = - 2$ and $D \cdot L = 1$); if there were a curve $F_{1}$ such
that
either $F_{1}^2  = 0, F_{1} \cdot (D - L) =1,2$ or
$F_{1}^2 = -2, F_{1} \cdot (D - L) = 0$ then $-4h-1$ would divide
$disc(F_{1},D-L) = -1, -4, -4h + 8$. Therefore $D - L$ is very ample
by
(i) of Lemma (3.5). If $h = 3$ we have $(D - L)^2 = 2$ and $D - L$
has no
base points and does not contract any curve by what we saw above,
hence we
can apply (iii) of Lemma (3.5). \ \qed \ (for Claim (3.6))

\noindent {\it Claim} (3.7). {\sl For $j=-1, k=2, h \geq 1$ we have
that
$D$ is very ample
for $h \geq 2$, $D$ defines a 2:1 finite morphism onto $\I^2$ for
$h=1$,
$D + L$ defines a birational morphism that contracts only the
irreducible
curve $L$.}

By Claim (3.7) we have that, for $j=-1, k=2, h \geq 1$, $H = A_{1} =
D +
D + L$ is very ample for $h \geq 2$ and also for $h=1$ by (i) of
Lemma
(3.5) since in this case $H^2=14$ and if there is a curve $F$ with
$F^2 = 0, F \cdot H =1,2$ or $F^2 = -2, F \cdot H = 0$ then
$disc \Gamma_{-1,2,1} = -8$ divides
$disc(F,H) =-1,-4, -28$. To see (ii) of Lemma (2.1) we check again
the
conditions
of Lemma (2.2): $(A_{1}+A_{2}+A_{3}) \cdot L = 2H \cdot L = 4$ and,
for
$h=1$,
$(A_{1}+A_{2}+A_{3}) \cdot A_{2} = 2H \cdot D = 12$.

\noindent {\it Proof of Claim (3.7):} If $h \geq 2$ we have $D^2
\geq 4$
hence $D$ is very ample since if there is a curve $F$ such that
$F^2 = 0, F \cdot D =1,2$ ($F \cdot D \not= 0$ since $D$ is ample)
then
$disc \Gamma_{-1,2,h}= -4h-4$ divides $disc(F,D) =-1, -4$. If $h=1$
and
there are irreducible curves $F_{1}$, $G$ and an integer $a \geq 2$
such
that $D \sim aF_{1} + G$, $F_{1}^2 = 0, G^2 = -2, F_{1} \cdot G =1$
then
$disc \Gamma_{-1,2,1}= -8$ divides $disc(F_{1},G) =-1$. Therefore
$D$ has
no base points and (iii) of Lemma (3.5) gives that it defines a 2:1
morphism onto $\I^2$ which is finite since $D$ is ample. Now notice
that
$L^2 =-2$ implies by Riemann-Roch that either
$L$ or $-L$ is effective, hence $L$ is because $D$ is nef and
$D \cdot L = 2$; also $L$ is irreducible for if $L = L_{1} + L_{2}$
then
necessarily $D \cdot L_{1} = D \cdot L_{2} = 1$ ($D$ is K\"ahler)
hence
$L_{1}$ and $L_{2}$ are irreducible and $-2 = L^2 = L_{1}^2 +
L_{2}^2 +
2L_{1} \cdot L_{2}$ gives that
$L_{1}^2 = L_{2}^2 = -2,  L_{1} \cdot L_{2} = 1$ ($L_{1} \not=
L_{2}$
since $L$ is indivisible) but then $disc \Gamma_{-1,2,h}= -4h-4$
divides
$disc(L_{1},L_{2}) = 3$, a contradiction. Since $L$ is irreducible
and
$(D + L) \cdot L = 0$ we have  that $D + L$ is nef and base point
free
because it has no base points on $L$ as it can be seen from the
exact
sequence $0 \to {\cal O}_{T_{-1,2,h}}(D) \to {\cal
O}_{T_{-1,2,h}}(D + L)
\to {\cal O}_{L} \to 0$ and the fact that $H^1({\cal
O}_{T_{-1,2,h}}(D)) =
0$ [SD, Proposition 2.6]. Also there is no curve $F_{2}$ with
$F_{2}^2 =
0, F_{2} \cdot (D+L) = 2$ else  $disc \Gamma_{-1,2,h}= -4h-4$
divides
$disc(F_{2},D+L) = -4$. By (iv) of Lemma (3.5) we get that
$D + L$
defines a birational morphism and if $F_{3}$ is an irreducible
contracted
curve  then $F_{3} \cdot (D+L) = F_{3} \cdot D + F_{3} \cdot L
> 0$
unless  $F_{3} = L$. \ \qed \ (for Claim (3.7))

\noindent {\it Claim} (3.8). {\sl For $j=0, k=1, h \geq 3$ we have
that
$D$ is very ample and $D - L + R$ defines a birational morphism
that
contracts only the irreducible curve $R$.}

By Claim (3.8) we have that, for $j=0, k=1, h \geq 3$, $H = A_{1} =
D + D - L + R$ is very ample and (ii) of Lemma (2.1) holds by Lemma
(2.2)
since $(A_{1}+A_{2}+A_{3}) \cdot R = 2H \cdot R = 4$.

\noindent {\it Proof of Claim (3.8):} First notice that $D$ is base
point
free, else there are irreducible curves $F$, $G$ and an integer $a
\geq 2$
such that $D \sim aF + G$, $F^2 = 0, G^2 = -2, F \cdot G =1$. But
then
$G \not= L$ (or $D - L$ is divisible) hence $1 = D \cdot L = aF
\cdot L
+ G \cdot L$ implies $F \cdot L = 0, G \cdot L = 1$ and then $disc
\Gamma_{0,1,h}= 8h + 10$ divides $disc(L,F,G) = 2$. To see that
$D$ is
very ample let $F_{1}$ be an irreducible curve such that $F_{1}^2
= 0,
F_{1} \cdot D = 2$ (the case $F_{1} \cdot D = 1$ is clearly
impossible).
Then $F_{1} \not= L$ and $0 \leq L \cdot F_{1} \leq 2$. In fact the
map
$\phi_{D}$ sends both $L$ and $F_{1}$ to lines in
$\I (H^0({\cal O}_{T_{0,1,h}}(D))^*)$. The possible values of
$disc(D,L,F_{1})$ are  $8, 12 - 2h, 16 - 8h$ and they are divisible
by
$disc \Gamma_{0,1,h}= 8h + 10$ only if $L \cdot F_{1} = 1, h = 6$.
Remark
that in this case $\phi_{D}$ is a finite morphism of $T_{0,1,6}$
onto a
rational normal scroll $V$ of degree $6$ in
$\I^7 = \I (H^0({\cal O}_{T_{0,1,6}}(D))^*)$, and $\phi_{D}(L)$ is
a
line directrix of $V$. One moment of reflection then shows that
$D = 5 F_{1} +
L + M$, where $M$ is a rational curve such that
$\phi_{D}(M) = \phi_{D}(L)$.
Since $1 = L \cdot D = 3 + L \cdot M$ we have $L = M$ and $D = 5
F_{1} +
2L$. But then $2 = D \cdot R = 5 F_{1} \cdot R + 2 L \cdot R = 5
F_{1}
\cdot R$, a contradiction. Hence the very ampleness of $D$ follows
by (i)
of Lemma (3.5). Note that $D$ very ample implies that $D - L$ is
nef.
Indeed the general curve in $|D - L|$ is irreducible by Bertini's
theorem
because $|D - L|$ has no base points and is connected since
$H^1({\cal
O}_{T_{-1,2,h}}(D - L)) = 0$ by [SD, Proposition 2.6]. Now notice
that $R$
is irreducible for if $R = R_{1} + R_{2}$ then $2 = D \cdot R$
implies $D
\cdot R_{1} = D \cdot R_{2} = 1$ and $R_{1}$, $R_{2}$ are
irreducible;
from $ R^2 = -2$ we get that $R_{1}^2 = R_{2}^2 = -2,  R_{1} \cdot
R_{2} =
1$ but then $disc  \Gamma_{0,1,h}= 8h+10$ divides $disc(D, R_{1},
R_{2}) =
6h+6$, which is impossible. Now $D - L$ nef implies that also $D -
L + R$
is nef since  $(D - L + R) \cdot R = 0$. Actually $|D - L + R|$
has no
base
points. This follows since $|D - L|$ has no base points and the
exact
sequence  $0 \to {\cal O}_{T_{-1,2,h}}(D - L) \to
{\cal O}_{T_{-1,2,h}}(D
- L + R) \to {\cal O}_{R} \to 0$ and
$H^1({\cal O}_{T_{-1,2,h}}(D - L)) =
0$ show that there are no base points on $R$ either. Now $(D - L
+ R)^2 =
2h - 2  \geq 4$ and suppose $F_{2}$ is an irreducible curve such
that
either $F_{2}^2 = -2, F_{2} \cdot (D - L + R) = 0$ or $F_{2}^2 =
0, F_{2}
\cdot  (D - L + R) = 2$. We will show that this is possible only
when
$F_{2} = R$, hence (iv) of Lemma (3.5) will give that $D - L + R$
defines
a birational morphism contracting only $R$. If $F_{2} \not= R$ in
the
first case we have $0 = F_{2} \cdot (D - L + R) = F_{2} \cdot (D
- L) +
F_{2} \cdot R$ hence $F_{2} \cdot (D - L) = F_{2} \cdot R = 0$ but
then
$disc \Gamma_{0,1,h}=  8h + 10$ divides $disc(D-L,R,F_{2}) = 8h-8$.
In the
second case we have  $2 = F_{2} \cdot (D - L) + F_{2} \cdot R$
hence
$F_{2} \cdot (D - L) = 0, 2, F_{2} \cdot R = 2, 0$ ($F_{2}$ is an
elliptic  curve hence we cannot have $F_{2} \cdot (D - L) = 1$) but
then
$disc \Gamma_{0,1,h}= 8h + 10$ divides $disc(D-L,R,F_{2}) = -8h+16,
8$,
both impossible.\ \qed \ (for Claim (3.8))

\noindent {\it Claim} (3.9). {\sl For $j=0, k=2, h \geq 1$ we have
that
$H = 2D + L + R$ is very
ample and $D + L$ (respectively $D + R$) defines a birational
morphism
that contracts only the irreducible curve $L$ (respectively $R$).}

By Claim (3.9) we see that, for $j=0, k=2, h \geq 1$, $H$ is very
ample
and (ii) of Lemma (2.1) holds by Lemma (2.2) since
$(A_{1}+A_{2}+A_{3})
\cdot L = (A_{1}+A_{2}+A_{3}) \cdot R = 2H \cdot L = 4$.

\noindent {\it Proof of Claim (3.9):} First we show that $L$ and
$R$ are
irreducible effective divisors. Suppose $L = L_{1} + L_{2}$ with
$L_{1},
L_{2}$ effective. Since $D \cdot L = 2$ we have $D \cdot L_{1} = D
\cdot
L_{2} = 1$ and $L_{1}$, $L_{2}$ are irreducible;  from $ L^2 = -2$
we get
that either $L_{1}^2 = -2,  L_{2}^2 = L_{1} \cdot L_{2} = 0$ or
$L_{1}^2 =
L_{2}^2 = -2,  L_{1} \cdot L_{2} = 1$.  But in these cases we have
that
$disc \Gamma_{0,2,h}= 8h+16$ divides $disc(D,  L_{1}, L_{2}) = 2,
6h+6$,
both impossible. Since $D \cdot L = D \cdot R$, $L^2 = R^2$ the same
proof
shows that $R$ is irreducible. To see that $H = 2D + L + R$ is very
ample
first notice that it is nef since $D$ is and $H \cdot L = H \cdot R
= 2$.
Let $F$ be an irreducible curve such that either $F^2 = - 2,
F \cdot H = 0$ or $F^2 = 0, F \cdot H = 1,2$. Then $F \not= L, R$
hence
$H \cdot F = 2D \cdot F + L \cdot F + R \cdot F \geq 2$, therefore
necessarily $D \cdot F = 1,  L \cdot F = R \cdot F = 0$ and $F^2 =
0$,
but then $disc \Gamma_{0,2,h} = 8h+16$ divides
$disc(D, L, F) = 2$. Therefore $H$ is very ample by (i) of Lemma
 (3.5).
Notice now that $D$ has no base points for otherwise we have
irreducible
curves $F_{1}$ and $G$ and an integer $a \geq 2$ such that $D
\sim
aF_{1} + G, F_{1}^2 = 0, G^2 = -2,  F_{1} \cdot G = 1$. But
then $L \not= F_{1}, G$ (else $D - L$ is divisible) and
$2 = D \cdot L = aF_{1} \cdot L + G \cdot L$ implies that either
$F_{1} \cdot L = 0, G \cdot L = 2$ or $F_{1} \cdot L = 1, G
\cdot L = 0$, but in both cases $disc(F_{1},L,G) = 2, 4$ is not
divisible by $disc \Gamma_{0,2,h}= 8h + 16$. $D$ being base point
free
yields that $D + L$ is also base point free. This follows from
$H^1({\cal O}_{T_{0,2,h}}(D)) = 0$ and the exact sequence
$0 \to {\cal O}_{T_{0,2,h}}(D) \to {\cal O}_{T_{0,2,h}}(D + L)
\to {\cal O}_{L} \to 0$. Since $(D + L)^2 = 2h + 2 \geq 4$ we use
now
(iv) of  Lemma (3.5) to conclude the proof of this Claim. Let
$F_{2}$
be an irreducible curve such that  either $F_{2}^2 = -2,  F_{2} \cdot
(D +
L) = 0$ or $F_{2}^2 = 0, F_{2} \cdot (D + L) = 2$. We will show that
this
is  possible only when $F_{2} = L$. If $F_{2} \not= L$ then $F_{2}
\cdot
(D + L) = F_{2} \cdot D + F_{2} \cdot L > 0$ hence $F_{2}^2 = 0$,
$F_{2}$
is an elliptic curve  and we cannot have $F_{2} \cdot D = 1$,
therefore we
have $F_{2} \cdot D = 2,  F_{2} \cdot L = 0$ but then $disc
\Gamma_{0,2,h}= 8h + 16$  divides $disc(D,L,F_{2}) = 8$. Again
replacing
$L$ by $R$ we get the statement for $D + R$. \ \qed \ (for Claim
(3.9))

\noindent {\it Claim} (3.10). {\sl For $j=1, k=5, h=2$ we have that
$D$ is
very ample and $L$ defines a 2:1 finite morphism onto $\I^2$.}

By Claim (3.10) we obtain that, for $j=1, k=5, h=2$, $H = A_{1} = D
+ L$
is very ample and (ii) of Lemma (2.1) holds by Lemma (2.2) since
$(A_{1}+A_{2}+A_{3}) \cdot A_{3} = 2H \cdot L = 2(D + L) \cdot L =
14.$

\noindent {\it Proof of Claim (3.10):} In this case the K3 surface
$T_{1,5,2}$ with Picard lattice $\Gamma_{1,5,2}$ can be taken to be
a
smooth quartic surface in $\I^3$ with Picard group generated by the
hyperplane section $D$ and a smooth irreducible genus $2$ quintic
curve
$L$ ($T_{1,5,2}$ exists by [Mo]). Moreover there is no irreducible
curve
$F$ such that $F^2 = -2,  F \cdot L = 0$ else $disc \Gamma_{1,5,2}=
-17$
divides $disc(L,F) = -4$. Therefore $L$ defines a 2:1 finite
morphism. \
\qed \ (for Claim (3.10))

\noindent {\it Claim} (3.11). {\sl For $j=1, k=4, h=1$ we have that
$H = D + L$ is very ample and both $D$ and $L$ define 2:1 finite
morphisms
onto $\I^2$.}

By Claim (3.11) we see that, for $j=1, k=4, h=1$, (ii) of Lemma
(2.1)
holds by Lemma (2.2) since $(A_{1}+A_{2}+A_{3}) \cdot A_{2} =
(A_{1}+A_{2}+A_{3}) \cdot A_{3} = 2H \cdot D = 2H \cdot L = 12.$

\noindent {\it Proof of Claim (3.11):} First of all notice that $D$
is
base point free otherwise there are irreducible curves $F, G$ and
an
integer $a \geq 2$ such that $D \sim aF + G$ but then $2 = D^2 = aF
\cdot
D + G \cdot D \geq 3$. Now let us see that $R$ is irreducible.
Suppose $R
= R_{1} + R_{2}$ with $R_{1}, R_{2}$ effective; since $D \cdot R
= 2$ we
have $D \cdot R_{1} = D \cdot R_{2} = 1$ hence $R_{1}$, $R_{2}$
are
irreducible and $R_{1}^2 = R_{2}^2 = - 2$ because $D$ has no base
points.
But then we have that $disc \Gamma_{1,4,1}= 30$ divides
$disc(D,R_{1},R_{2}) = 12$, a contradiction. Also $L$ can be
assumed
irreducible. To see this observe first that $R$ is not a fixed
component
of the linear system $|L|$. Indeed $L - R$ is effective because
$(L - R)^2 = - 2, (L - R) \cdot D = 2$ and connected since if we
had
$L - R = B_{1} + B_{2}$ with $B_{1}, B_{2}$ effective then
$D \cdot B_{1} = D \cdot B_{2} = 1$ hence $B_{1}$, $B_{2}$ are
irreducible
and distinct (because $L - R$ is indivisible), therefore $B_{1}
\cdot
B_{2} \geq 0$ and if $B_{1} \cdot B_{2} = 0$, since
$B_{1}^2 = B_{2}^2 = - 2$, we get that $disc \Gamma_{1,4,1}= 30$
divides
$disc(D,B_{1},B_{2}) = 12$, a contradiction. Since
$h^0({\cal O}_{T_{1,4,1}}(L)) \geq 3,
h^0({\cal O}_{T_{1,4,1}}(L-R)) = 1$
we deduce that $R$ is not a fixed component of $|L|$. Let
$L' = L_{1} + L_{2}$ be a divisor of $|L|$ not containing $R$,
with $L_{1},
L_{2}$ effective. We have $R \cdot L_{i} \geq 0, i = 1, 2,$ hence
$0
\leq R \cdot L_{i} \leq 1$ and $1 \leq D \cdot L_{i} \leq 3$. If
$D \cdot
L_{1} = 1$ then $disc \Gamma_{1,4,1}= 30$ divides
$disc(D,L_{1},R) = 18, 20$ which is impossible. If $D \cdot L_{1}
= 3$ then
$D
\cdot L_{2} = 1$ and we conclude as above replacing $L_{1}$ with
$L_{2}$.
Hence $D \cdot L_{1} = D \cdot L_{2} = 2$ and they are both
irreducible.
The linear system cut out on $L_{1}$ by $|D|$ has dimension $2$
otherwise
$D - L_{1}$ is effective and $D \cdot (D - L_{1}) = 0$ hence $D =
L_{1}$
and $-1 = R \cdot (L - D) = R \cdot L_{2} \geq 0$, a contradiction.
This
implies that $L_{1}$ is rational and therefore $L_{1}^2 = - 2$ and
the
same holds for $L_{2}$. Therefore we can assume $R \cdot L_{1} = 0$
and
this is impossible because $disc \Gamma_{1,4,1}= 30$ does not divide
$disc(D,L_{1},R) = 24$. This shows that $L$ is irreducible.
Since  $L^2 = 2$ we have that $L$ is nef and so is $H = D + L$. Let
$F$ be
an irreducible curve such  that either $F^2 = -2, F \cdot H = 0$ or
$F^2 =
0, F \cdot H = 1,2$.  Then $F \not= L$ hence $H \cdot F = D \cdot F
+ L
\cdot F \geq 1$ and therefore we have $F^2 = 0$ and either $D \cdot
F = 1,
L \cdot F = 0$ or  $D \cdot F =  L \cdot F = 1$ or $D \cdot F = 2, L
\cdot
F = 0$, but then $disc \Gamma_{1,4,1}= 30$ divides $disc(D,F,L) =
 -2, 4,
-8$, a contradiction. Therefore $H$ is very ample by (i) of Lemma
(3.5).
Since $D$ has no base points, by (ii) and
(iii) of Lemma (3.5) we have that $D$ defines a 2:1 finite morphism
onto
$\I^2$ because $D$ is ample. Similarly for $L$ (which is irreducible,
hence
base point free) there is no  irreducible curve $F_{2}$ such that
$F_{2}
\cdot L  = 0, F_{2}^2 = -2$ else setting $x = R \cdot F_{2}$ we get
that
$disc \Gamma_{1,4,1}= 30$ divides $disc(R,L,F_{2}) = -2x^2 + 10$, but
then
$2x^2 \equiv 4 \ (mod \ 6)$ and this is not possible for an integer
$x$.
\ \qed \ (for Claim (3.11))

\noindent {\it Claim} (3.12). {\sl For $j=1,2, k \geq j+4, h=2,3$ we
have that $D$ is very ample and so is $L$, except for $j=1$, and in
that
case $L$ defines a 2:1 finite morphism onto $\I^2$.}

 From Claim (3.12) we see that, for $j=1,2, k \geq j+4, h=2,3$, $2D+L$
and
$D + 2L$ are very ample and both satisfy (ii) of Lemma (2.1) by Lemma
(2.2) since $(A_{1}+A_{2}+A_{3}) \cdot L = (2D+L) \cdot L = 2k + 2j
\geq
4j+8 \geq 12$ or  $(A_{1}+A_{2}+A_{3}) \cdot L = (D+2L) \cdot L =
k+4j
\geq 5j+4 \geq 9$.

\noindent {\it Proof of Claim (3.12):} If there is a curve $F$ such
that
$F^2 = 0, F \cdot D =1,2$ ($F \cdot D = 0$ is excluded since
$D$ is an ample class) then $disc \Gamma_{jkh} = 8j-k^2$ divides
$disc(F,D) =-1,-4$, but $k^2 - 8j \geq 13$ does not divide $1$ nor
$4$.
Hence $D$ is very ample by (i) of Lemma (3.5). When $h = 2$ it is
shown
in [Mo] that $L$ can be assumed to be irreducible base point free.
Hence,
by (iii) of Lemma (3.5), it defines a 2:1 finite morphism onto
$\I^2$ for
$j=1$ since there  is no irreducible curve $F_{1}$ such that $F_{1}
\cdot
L = 0,  F_{1}^2 = -2$ else  $disc \Gamma_{1,k,2}= 8 - k^2 \leq -25$
divides $disc(L,F_{1}) = -4$. For $j=2$ if there were a curve
$F_{2}$ such
that either $F_{2}^2 = 0, F_{2} \cdot L =1,2$ or  $F_{2}^2 = -2,
F_{2}
\cdot L = 0$ then $disc(F_{2},L) =-1, -4, -8$ would  be divisible by
$disc
\Gamma_{2,k,2} = 16 - k^2$. But $k^2 - 16 \geq 20$, therefore $L$ is
very ample by (i) of Lemma (3.5). For $h = 3$ it suffices to show
that
$L$ is irreducible: In fact it is then nef, base  point free and
defines a
2:1 finite morphism onto $\I^2$ again by (iii)  of Lemma (3.5) since
if
there were a curve $F_{3}$ such that $F_{3}^2 =  -2, F_{3} \cdot L =
0$
then $disc(F_{3},L) = -4$ would be  divisible by $disc
\Gamma_{1,5,3} =
-13$. To show that $L$ can be assumed to be irreducible notice that
$D -
L$ is effective by Riemann-Roch since  $(D - L)^2 = - 2, (D - L)
\cdot D =
1$ and in fact it is a line $L'$ in  the embedding given by $D$.
Now a
general element $L \in |D - L'|$ is  smooth since $D - L'$ is base
point
free and connected (as in the proof of Claim (3.8)). \ \qed \ (for
Claim
(3.12))

The proof of Proposition (3.1) is now complete. \ \qed

\vskip .5cm
{\bf REFERENCES}
\baselineskip 12pt
\vskip .5cm
\item{[BEL]} Bertram,A., Ein,L., Lazarsfeld,R.:\ Surjectivity of
Gaussian
maps for line bundles of large degree on curves.\ In:\ {\it Algebraic
Geometry, Proceedings Chicago 1989.\ Lecture Notes in Math.\ \bf
1479}.\
Springer, Berlin-New York:\ 1991, 15-25.
\vskip .2cm
\item{[BPV]} Barth,W., Peters,C., Van de Ven,A.:\ Compact complex
surfaces.\ {\it Ergebnisse der Mathematik \bf 4}.\ Springer,
Berlin-New
York:\ 1984.
\vskip .2cm
\item{[CLM]} Ciliberto,C., Lopez,A.F., Miranda,R.:\ Projective
degenerations of K3 surfaces, Gaussian maps and Fano threefolds.\
{\it
Invent.\ Math.\ \bf 114}, (1993) 641-667.
\vskip .2cm
\item{[D]} Duflot,J.:\ Gaussian maps for double coverings.\
{\it Manuscripta Math.\ \bf 82}, (1994) 71-87.
\vskip .2cm
\item{[EEK]} Ein,L., Eisenbud,D., Katz,S.:\ Varieties cut out by
quadrics:\ Scheme-theoretic versus homogeneous generation of ideals.\
In:\ {\it Algebraic Geometry, Proceedings Sundance 1986.\ Lecture
Notes
in Math.\ \bf 1311}.\ Springer, Berlin-New York:\ 1988, 51-70.
\vskip .2cm
\item{[G]} Griffiths,P.:\ Hermitian differential geometry, Chern
classes,
and positive vector bundles.\ In:\ {\it Global Analysis, Papers in
honor
of K.\ Kodaira}.\ Princeton University Press:\ Princeton, NJ, 1969,
181-251.
\vskip .2cm
\item{[K]} Kumar,S.:\  Proof of Wahl's conjecture on surjectivity
of the Gaussian map for flag varieties.\ {\it Amer.\ J.\ Math.\
\bf 114},
(1992) 1201-1220.
\vskip .2cm
\item{[Ma]} Mayer,A.:\ Families of K3 surfaces.\ {\it Nagoya Math.\
J.\
\bf 48}, (1972) 1-17.
\vskip .2cm
\item{[Mo]} Mori,S.:\ On degrees and genera of curves on smooth
quartic
surfaces in $\I^3$.\ {\it Nagoya Math.\ J.\ \bf 96}, (1984) 127-132.
\vskip .2cm
\item{[SD]} Saint-Donat,B.:\ Projective models of K3 surfaces.\
{\it
Amer.\ J.\ Math.\ \bf 96}, (1974) 602-639.
\vskip .2cm
\item{[W1]} Wahl,J.:\ The Jacobian algebra of a graded Gorenstein
singularity.\ {\it Duke Math.\ J.\ \bf 55}, (1987) 843-871.
\vskip .2cm
\item{[W2]} Wahl,J.:\ Gaussian maps and tensor products of
irreducible
representations.\ {\it Manuscri-}
\noindent {\it pta Math.\ \bf 73}, (1991) 229-259.
\vskip .2cm
\item{[W3]} Wahl,J.:\ Introduction to Gaussian maps on an algebraic
curve.\ In: {\it Complex Projective Geometry, Trieste-Bergen 1989.
London
Math.\ Soc.\ Lecture Notes Series \ \bf 179}. Cambridge Univ.\
Press:\
1992, 304-323.
\vskip .2cm
\item{[Z]} Zak,F.L.:\ Some properties of dual varieties and their
application in projective geometry.\ In:\ {\it Algebraic Geometry,
Proceedings Chicago 1989.
\ Lecture Notes in Math.\ \bf 1479}.\ Springer, Berlin-New York:\
1991,
273-280.
\vskip .5cm
ADDRESSES OF THE AUTHORS:
\vskip .2cm
Ciro Ciliberto, Dipartimento di Matematica, Universit\`a di Roma II,
Tor Vergata,

Viale della Ricerca Scientifica, 00133 Roma, Italy

e-mail: ciliberto@mat.utovrm.it
\vskip .2cm
Angelo Felice Lopez, Dipartimento di Matematica, Terza Universit\`a
di Roma,

Via Corrado Segre 2, 00146 ROMA, Italy

e-mail: lopez@matrm3.mat.uniroma3.it
\vskip .2cm
Rick Miranda, Department of Mathematics, Colorado State University,
Ft. Collins,

CO 80523, USA

e-mail: miranda@math.colostate.edu

\end